\documentstyle[11pt]{article}

\begin{document}
\pagestyle{plain} \headheight=5mm \footheight=5mm \topmargin=-5mm

\title{A Poincar\'e-Hopf type formula for\\ Chern character numbers}
\author{Huitao Feng\thanks{Partially supported by NNSFC, MOEC and NSFC.},\quad Weiping Li\quad  and \quad Weiping Zhang\thanks{Partially supported by NNSFC and MOEC.}}
\date{}

\maketitle

\newtheorem{Def}{Definition}[section]
\newtheorem{Th}{Theorem}[section]
\newtheorem{Prop}{Proposition}[section]
\newtheorem{Not}{Notation}[section]
\newtheorem{Lemma}{Lemma}[section]
\newtheorem{Rem}{Remark}[section]
\newtheorem{Cor}{Corollary}[section]

\def\s{\section}
\def\ss{\subsection}

\def\d{\begin{Def}}
\def\t{\begin{Th}}
\def\p{\begin{Prop}}
\def\n{\begin{Not}}
\def\la{\begin{Lemma}}
\def\r{\begin{Rem}}
\def\c{\begin{Cor}}
\def\ee{\begin{equation}}
\def\aa{\begin{eqnarray}}
\def\ya{\begin{eqnarray*}}
\def\bd{\begin{description}}

\def\ed{\end{Def}}
\def\et{\end{Th}}
\def\epo{\end{Prop}}
\def\en{\end{Not}}
\def\el{\end{Lemma}}
\def\er{\end{Rem}}
\def\ec{\end{Cor}}
\def\eee{\end{equation}}
\def\eaa{\end{eqnarray}}
\def\ey{\end{eqnarray*}}
\def\ebd{\end{description}}

\def\nn{\nonumber}
\def\bp{{\bf Proof.}\hspace{2mm}}
\def\qe{\hfill$\Box$}
\def\lj{\langle}
\def\rj{\rangle}
\def\dd{\diamond}
\def\ox{\mbox{}}
\def\lb{\label}
\def\rel{\;{\rm rel.}\;}
\def\vp{\varepsilon}
\def\ep{\epsilon}
\def\mod{\;{\rm mod}\;}
\def\exp{{\rm exp}\;}
\def\Lie{{\rm Lie}}
\def\dim{{\rm dim}}
\def\im{{\rm im}\;}
\def\Lag{{\rm Lag}}
\def\Gr{{\rm Gr}}
\def\Sym{{\rm Sym}}
\def\span{{\rm span}}
\def\Spin{{\rm Spin}}
\def\sign{{\rm sign}\;}
\def\Supp{{\rm Supp}\;}
\def\Sp{{\rm Sp}\;}
\def\ind{{\rm ind}\;}
\def\rank{{\rm rank}\;}
\def\Sg{{\Sp(2n,\C)}}
\def\Na{{\cal N}}
\def\det{{\rm det}\;}
\def\dist{{\rm dist}}
\def\deg{{\rm deg}}
\def\tr{{\rm tr}\;}
\def\tr_s{{\rm tr}_s\;}
\def\ker{{\rm ker}\;}
\def\Vect{{\rm Vect}}
\def\H{{\bf H}}
\def\K{{\rm K}}
\def\R{{\bf R}}
\def\C{{\bf C}}
\def\Z{{\bf Z}}
\def\N{{\bf N}}
\def\F{{\bf F}}
\def\Da{{\bf D}}
\def\A{{\bf A}}
\def\La{{\bf L}}
\def\x{{\bf x}}
\def\y{{\bf y}}
\def\Ga{{\cal G}}
\def\Ha{{\cal H}}
\def\L{{\cal L}}
\def\Pa{{\cal P}}
\def\Ua{{\cal U}}
\def\E{{\rm E}}
\def\J{{\cal J}}

\def\m{{\rm m}}
\def\ch{{\rm ch}}
\def\gl{{\rm gl}}
\def\Gl{{\rm Gl}}
\def\Sp{{\rm Sp}}
\def\sf{{\rm sf}}
\def\U{{\rm U}}
\def\O{{\rm O}}
\def\F{{\rm F}}
\def\P{{\rm P}}
\def\D{{\rm D}}
\def\T{{\rm T}}
\def\Sa{{\rm S}}

\begin{abstract} For two complex  vector bundles admitting a homomorphism with isolated singularities between them,
 we establish  a Poincar\'e-Hopf type formula for the difference of the  Chern character numbers of these two vector bundles. As a consequence, we extend the
 original Poincar\'e-Hopf index  formula to the case of complex vector fields.\end{abstract}

\s{Introduction and the statement of the main result}

\quad\ Let $M$ be a closed, oriented, smooth manifold of dimension
$2n$. Let $E_+$, $ E_-$ be two complex vector bundles over $M$.

Let $v\in\Gamma({\rm Hom}(E_+,E_-))$ be a homomorphism between $E_+$
and $E_-$. Let $Z(v)$ denote the set of the points at which $v$ is
singular (that is, not invertible). We assume that the following
basic assumption holds.

$\ $

\noindent{\bf Basic Assumption 1.0.}  The point set $Z(v)$ consists
of a finite number of points in $M$.

\bigskip
  For any $p\in Z(v)$,
we choose a small open ball $B(p)$ centered at $p$ such that the
closure $\overline{B(p)}$ contains no points in $Z(v)\setminus p$.
Then, when restricted  to the boundary $\partial{B(p)}$, the linear
map
$$v|_{\partial{B(p)}}:{E_+}|_{\partial{B(p)}}\to
{E_-}|_{\partial{B(p)}},\eqno(1.1)$$ which we denote by $v_p$, is
invertible. The map $v_p$ determines an element in $K^1(S^{2n-1})=
{\bf Z}$ which we denote by ${\rm deg}(v_p)\in {\bf
Z}$.\footnote{One way to define ${\rm deg}(v_p)$ is that $v_p$ in
(1.1) defines a complex vector bundle $E_{v(p)}$ over a sphere
$S^{2n}(v_p)$ with $\partial{B(p)}$ as an equator. Then one can
define ${\rm deg}(v_p)=\langle {\rm
ch}(E_{v(p)}),[S^{2n}(v_p)]\rangle$.}

The main result in this paper is the following theorem:

\bigskip
\noindent {\bf Theorem 1.1.} Under the Basic Assumption 1.0,  the
following identity holds,
$$\left\lj{\rm ch}(E+)-{\rm ch}(E_-),[M]\right\rj=(-1)^{n-1}\sum_{p\in Z(v)}{\rm deg}\left(v_p\right).
\eqno(1.2)$$

Our original motivation is to establish an extension of the
Poincar\'e-Hopf index formula for vector fields with isolated zero
points (cf. \cite[Theorem 11.25]{BT}) to the case of complex
vector fields, under the framework considered by Jacobowitz in
\cite{J}.

To be more precise, let $T_{\bf C}M=TM\otimes {\bf C}$ denote the
complexification of the tangent vector bundle $TM$. Let
$K=\xi+\sqrt{-1}\eta\in\Gamma(T_{\bf C}M)$ be a smooth section of
$T_{\bf C}M$, with $\xi,\,\eta\in\Gamma(TM)$.

Let $g^{TM}$ be a Riemannian metric on $TM$, then it induces
canonically a complex symmetric bilinear form $h^{T_{\bf C}M}$ on
$T_{\bf C}M$, such that
$$h^{T_{\bf
C}M}\left(K,K\right)=|\xi|^2_{g^{TM}}-|\eta|^2_{g^{TM}}+2\sqrt{-1}\langle
\xi,\eta\rangle_{g^{TM}}.\eqno (1.3)$$

Jacobowitz proved in \cite{J} the following vanishing result.

$\ $

\noindent{\bf Proposition 1.2.} (Jacobowitz  \cite{J}) If $h^{T_{\bf
C}M}(K,K )$ is nowhere zero on $M$, then the Euler number of $M$
vanishes: $\chi(M)=0$.

$\ $

If one takes $\eta=0$, then Proposition 1.2 reduces to the classical
Hopf vanishing result: $\chi(M)=0$ if $M$ admits a nowhere zero
vector field.

Jacobowitz asked in \cite{J} whether there is a counting formula for
$\chi(M)$ of Poincar\'e-Hopf type, extending Proposition 1.2 to the
case where $h^{T_{\bf C}M}(K,K )$ vanishes somewhere on $M$. In
Section 3, we will establish such a formula as an application of
Theorem 1.1, while Theorem 1.1 itself will be proved in Section 2.

\s{A Proof of Theorem 1.1}

We will use the superconnection formalism developed in \cite{Q} to
prove Theorem 1.1.

Due to  the topological nature of both sides of (1.2), we first make
some simplifying assumptions on the metrics and connections near the
set of singularities $Z(v)$.

First of all, we assume that there is a Riemannian metric $g^{TM}$
on $TM$ such that for any $p\in Z(v)$, there is a coordinate system
$\left(x_1,\cdots, x_{2n}\right)$, with $0\leq x_i\leq 1$ for $1\leq
i\leq 2n$,  centered around $p$ such that
$$B_p (1)=\left\{(x_1,\ldots,x_{2n})|\sum^{2n}_{i=1}x^2_i\leq 1 \right\}
\subset M\setminus(Z(v)\setminus\{p\}) \eqno(2.1)$$ and
$$\left. g^{TM}\right|_{B_p(1)}=dx_1^2+dx_2^2+\cdots+dx_{2n}^2,
\eqno(2.2)$$ that is, the metric $g^{TM}$ is Euclidean on each
$B_p(1)$, $p\in Z(v)$.

On the other hand, on each $B_p(1)$, the bundles $E_\pm$ are trivial
vector bundles.  We equip these two trivial vector bundles over
$B_p(1)$ the trivial metrics and trivial connections respectively.
Moreover, we can deform $v$ near $\partial B_p (1)$, so that
$v_p:E_+|_{\partial B_p (1)}\rightarrow E_-|_{\partial B_p (1)}$ is
{\it unitary}, while still keep the   new homomorphism nonsingular
on $M\setminus Z(v)$.

By partition of unity, we may then construct Hermitian metrics and
connections $\nabla^{E_\pm}$ on $E_\pm$ over $M$ such that the above
simplifying assumptions hold on $\cup_{p\in Z(v)}B_p(1)$.

We now follow the formalism in \cite{Q}.

Let $E=E_+\oplus E_-$ be the ${\bf Z}_2$-graded complex vector
bundle over $M$. Let $\nabla^E=\nabla^{E_+}\oplus \nabla^{E_-}$ be
the ${\bf Z}_2$-graded connection on $E$.

Let $v:E_+\rightarrow E_-$  extend to an (odd) endomorphism of $E$
by acting as zero on $E_-$, with the notation unchanged. Let
$v^*:E_-\rightarrow E_+$ (and thus also extends to an (odd)
endomorphism of  $E$) be the adjoint of $v$ with respect to the
Hermitian metrics on $E_\pm$ respectively.

Set $V=v+v^*$. Then $V$ is an odd endomorphism of $E$. Moreover,
$V^2$ is fiberwise positive over $M\setminus Z(v)$.

We fix a square root of $\sqrt{-1}$. Let
$\varphi:\Omega^*(M)\rightarrow \Omega^*(M)$ be the rescaling  on
differential forms such that for any differential form $\alpha$ of
degree $k$, $\varphi(\alpha)=(  2\pi\sqrt{-1})^{-{k\over 2}}\alpha$.
The final formulas below will not depend on the choice of this
square root.

 For any $t\in{\bf R}$, let ${\bf A}_t$ be the
superconnection on $E$, in the sense of Quillen \cite{Q}, defined by
$${\bf A}_t=\nabla^E+tV.
\eqno(2.3)$$  Let ${\rm ch}(E,{\bf A}_t)$ be the associated Chern
character form defined by
$${\rm ch}\left(E,{\bf A}_t\right)=\varphi\,\tr_s\left[e^{-{\bf
A}_t^2}\right].\eqno(2.4)$$ The following transgression formula has
been proved in \cite[(2)]{Q},
$$
{{\partial{\rm ch}\left(E,{\bf A}_t\right)}\over{\partial t}}
=-{{1\over{{\sqrt{2\pi{\sqrt{-1}}}}}}}d\,
\varphi\,\tr_s\left[Ve^{-{\bf A}_t^2}\right]. \eqno(2.5)$$

 Set for
any $T>0$,
$$\gamma(T)={{1\over{{\sqrt{2\pi{\sqrt{-1}}}}}}}\varphi\int^T_0\tr_s\left[Ve^{-{\bf A}_t^2}\right]dt.
\eqno(2.6)$$

 From (2.5) and (2.6), one gets
$${\rm ch}\left(E,{\bf A}_0\right)-{\rm ch}\left(E,{\bf A}_T\right)=d\gamma(T).
\eqno(2.7)$$

Set ${  M}_1=M\setminus\bigcup_{p\in Z(v)}B_{p}(1).  $

Since $V$ is invertible on $M_1$, by proceeding as in \cite[\S
4]{Q}, one sees that the following identity holds uniformly on
$M_1$,
$$\lim_{T\to+\infty}{\rm ch}\left(E,{\bf A}_T\right)=0.
\eqno(2.8)$$

\noindent{\bf Lemma 2.1.} The following identity holds,
$$\left\lj{\rm ch}\left(E_+\right)-{\rm
ch}\left(E_-\right),[M]\right\rj=-\sum_{p\in Z(v)}\lim_{T\to
+\infty}\int_{\partial{B_p(1)}}\gamma(T). \eqno(2.9)$$ \bp Since by
our choice the   connections $\nabla^{E_\pm}$ are the trivial
connections when restricted to $\bigcup_{p\in Z(v)}B_p(1)$,   one
has
$$\left\lj{\rm ch}\left(E_+\right)-{\rm
ch}\left(E_-\right),[M]\right\rj=\int_M{\rm ch}\left(E,{\bf
A}_0\right)=\varphi\int_M\tr_s\left[e^{-(\nabla^E)^2}\right]
=\varphi\int_{{ M}_1}\tr_s\left[e^{-\left(\nabla^E\right)^2}\right].
\eqno(2.10)$$ By (2.7), (2.8) and (2.10), we have
\begin{eqnarray*}
\left\lj{\rm ch}\left(E_+\right)-{\rm
ch}\left(E_-\right),[M]\right\rj%&=&\int_{{ M}_1}{\rm ch}\left(E,{\bf A}_0\right)\\
&=&\lim_{T\to+\infty}\left(\int_{{ M}_1}{\rm ch}\left(E,{\bf
A}_0\right)
-\int_{{  M}_1}{\rm ch}\left(E,{\bf A}_T\right)\right)\\
&=&\lim_{T\to+\infty}\int_{{  M}_1}d\gamma(T)
=\lim_{T\to+\infty}\int_{\partial{  M}_1}\gamma(T)\\
&=&-\sum_{p\in Z(v)}\lim_{T\to
+\infty}\int_{\partial{B_{p}(1)}}\gamma(T),
\end{eqnarray*} where the last equality comes from the orientation
consideration. \ \ Q.E.D.

$\ $

 Recall that the map $v_p$ is the restriction of $v$ on
$\partial B_{p}(1)$ (cf. (1.1)).

\bigskip

\noindent {\bf Lemma 2.2.} For any $p\in Z_v$, the following
identity holds, $$\lim_{T\to
+\infty}\int_{\partial{B_{p}(1)}}\gamma(T)= (-1)^{n }\deg(v_p).
\eqno(2.11)$$

\noindent\bp For any $p\in Z(v)$, since when restricted on the
 sphere $\partial{B_{p}(1)}$, the homomorphism  $v$ has been deformed to be unitary,
we get that $v^*=v^{-1}$ and $V^2$ is the identity map acting on
$E|_{\partial{B_{p}(1)}}$. Also, since $\nabla^E$ is the trivial
connection over $B_{p}(1)$, we will use the simplified notation $d$
for it. By (2.3), one has on $B_p(1)$ that
$$A_t=d+tV,\ \ \ A_t^2=d^2+t[d,V]+t^2V^2=t^2{\rm Id}_E+tdV.$$

 One  then deduces that
\begin{eqnarray*}
&&\int_{\partial{B_{p}(1)}}\gamma(T)
={{1\over{{\sqrt{2\pi{\sqrt{-1}}}}}}}\varphi\int_{\partial{B_{p}(1)}}\int^T_0\tr_s\left[Ve^{-{\bf A}_t^2}\right]dt\\
&=&{{1\over{{\sqrt{2\pi{\sqrt{-1}}}}}}}\varphi\int_{\partial{B_{p}(1)}}\int^T_0e^{-t^2}
\tr_s\left[Ve^{-tdV}\right]dt\\
&=&{1\over{(2\pi{\sqrt{-1}})^n}}{{-1}\over{(2n-1)!}}\int^T_0t^{2n-1}e^{-t^2}dt
\int_{\partial{B_{p}(1)}}\left({\rm tr}_{E_+}\left[v^*dv\left(dv^*dv\right)^{n-1}\right]-{\rm tr}_{E_-}\left[vdv^*\left(dvdv^*\right)^{n-1}\right]\right)\\
%&=&{1\over{(2\pi{\sqrt{-1}})^n}}{{(-1)^n}\over{(2n-1)!}}\int^T_0t^{2n-1}e^{-t^2}dt
% \cdot
%\int_{\partial{B_{p}(1)}}\left({\rm tr}\left[v^{-1}dv\left(v^{-1}dv\right)^{2n-2}\right]+{\rm tr}\left[(dv)v^{-1}\left((dv)v^{-1}\right)^{2n-2}\right]\right)\\
&=&{1\over{(2\pi{\sqrt{-1}})^n}}{{2(-1)^n}\over{(2n-1)!}}\int^T_0t^{2n-1}e^{-t^2}dt
\int_{\partial{B_{p}(1)}}{\rm
tr}_{E_+}\left[\left(v^{-1}dv\right)^{2n-1}\right].
\end{eqnarray*}
%where the trace is taken on $E_+|_{\partial B_p(1)}$.

 Hence,
\begin{eqnarray*}
\lim_{T\to +\infty}\int_{\partial{B_{\epsilon}(p)}}\gamma(T)
&=&{1\over{(2\pi{\sqrt{-1}})^n}}{{2(-1)^n}\over{(2n-1)!}}\int^{+\infty}_0t^{2n-1}e^{-t^2}dt
\int_{\partial{B_{p}(1)}}{\rm tr}_{E_+}\left[\left(v^{-1}dv\right)^{2n-1}\right]\\
&=&{1\over{(2\pi{\sqrt{-1}})^n}}{{(-1)^n(n-1)!}\over{(2n-1)!}}\int_{\partial{B_{p}(1)}}{\rm tr}_{E_+}\left[\left(v^{-1}dv\right)^{2n-1}\right]\\
%&=&{{(-1)^{n-1}}\over{(-2\pi{\sqrt{-1}})^n}}\int_{\partial{B_{p}(1)}}\left((-1)^{n-1}{{(n-1)!}\over{(2n-1)!}}
%{\rm tr}\left[\left(v^{-1}dv\right)^{2n-1}\right]\right)\\
&=& (-1)^{n}\deg(v_p),
\end{eqnarray*}
where one compares with \cite[Propositions 1.2 and 1.4]{G} for the
last equality.\ \ Q.E.D.

\bigskip
  From  Lemmas
2.1 and 2.2, one gets Theorem 1.1. \ \ Q.E.D.

$\ $

We conclude this section with the following result which is
complementary to Theorem 1.1.

$\ $

\noindent {\bf Lemma 2.3.} Under the Basic Assumption 1.0, for any
   closed form $\alpha\in \Omega^*(M)$ without degree zero component, one has
   $$\langle [\alpha]\left({\rm ch}\left(E_+\right)-{\rm ch}\left(E_-\right)\right),[M]\rangle =0,\eqno(2.12)$$
   where $[\alpha]\in H^*(M,{\bf C})$ is the de Rham cohomology
   class induced by $\alpha$.

 \noindent{\bf Proof.} By the Poincar\'e lemma (cf. \cite{BT}), as $\alpha$ is
   closed and contains no zero degree component, on each $B_p(1)$,  $p\in Z(v)$,
   there exists a form $\beta_p$ such that $\alpha=d\beta_p$ on an
   open neighborhood of
   $B_p(1)$.

   By partition of unity, one then constructs a differential form
   $\beta$ on $M$ such that $\beta=\beta_p$ on   each $B_p(1)$, $p\in Z(v)$. Then,
   $$\alpha-d\beta =0\eqno(2.13)$$ on  $\cup_{p\in Z(v)}B_p(1)=M\setminus M_1$.

   On the other hand, by (2.4) and (2.5) one knows that for any $t\geq 0$, one
   has
   $$\langle [\alpha]\left({\rm ch}\left(E_+\right)-{\rm ch}\left(E_-\right)\right),[M]\rangle =\int_M(\alpha-d\beta)\varphi\, {\rm
   tr}_s\left[e^{-{\bf A}_t^2}\right].\eqno(2.14)$$

   From (2.8), (2.13) and (2.14), and by taking $t\rightarrow +\infty$, one gets (2.12).\ \ Q.E.D.

\s{A Poincar\'e-Hopf formula for complex vector fields}

Let $M$ be a closed and oriented manifold of dimension $2n$. Let
$g^{TM}$ be a Riemannian metric on $TM$. Let $T_{\bf C}M=TM\otimes
{\bf C}$ be the complexification of $TM$. Then $g^{TM}$ extends to a
symmetric bilinear form $h^{T_{\bf C}M}$ on $T_{\bf C}M$.

Let $K=\xi+\sqrt{-1}\eta\in \Gamma(T_{\bf C}M)$ be a complex  vector
field on $M$, with $\xi,\,\eta\in\Gamma(TM)$. Then one has
$$h^{T_{\bf
C}M}(K,K)=|\xi|_{g^{TM}}^2-|\eta|_{g^{TM}}^2+\sqrt{-1}\langle\xi,\eta\rangle_{g^{TM}}.\eqno(3.1)$$

Let $Z_K$ be the zero set of $h^{T_{\bf C}M}(K,K)$, that is,
$$Z_K=\left\{x\in M:h^{T_{\bf
C}M}(K(x),K(x))=0\right\}.\eqno(3.2)$$

In the rest of this section, we make the following assumption.

$\ $

\noindent{\bf Basic Assumption 3.0.} The set $Z_K$ consists of a
finite number of points.

$\ $

Let $a_0>0$ be the  injectivity radius  of $g^{TM}$. Let
$0<\epsilon<{a_0\over 2}$.

For any $p\in Z_K$, let $B_p(\epsilon)=\{x\in M: d^{g^{TM}}(x,p)\leq
\epsilon\}$ be the Riemannian ball centered at $p$. We may take
$\epsilon$ small enough so that each $B_p(\epsilon)$ does not
contain points in $Z_K\setminus \{ p\}$.

Let $S(TB_p(\epsilon))=S_+(TB_p(\epsilon))\oplus
S_-(TB_p(\epsilon))$ be the Hermitian  bundle of spinors associated
with $(TB_p(\epsilon),g^{TM}|_{B_p(\epsilon)})$. Let $\tau$ be the
the involution on $S(TB_p(\epsilon))$ such that
$\tau|_{S_\pm(TB_p(\epsilon))}=\pm{\rm
Id}|_{S_\pm(TB_p(\epsilon))}$. Let $c(\cdot)$ denote the Clifford
action on $S(TB_p(\epsilon))$.\footnote{For a thorough treatment of
 spin geometry involved here, see \cite{LM}.}

Let $v_K(p): \Gamma(S_+(TB_p(\epsilon)))\rightarrow
\Gamma(S_-(TB_p(\epsilon)))$ be defined by $$v_K(p)=\tau
c(\xi)+\sqrt{-1}c(\eta).\eqno(3.3)$$ Then one can prove (see Lemma
3.2 below) that the restriction of $v_K(p)$ on the sphere $\partial
B_p(\epsilon)$ is invertible. Thus it defines an integer ${\rm
deg}(v_K(p)|_{\partial B_p(\epsilon)})\in {\bf Z}=K^1(\partial
B_p(\epsilon))$.

We can now state the main result of this section as follows.

$\ $

\noindent {\bf Theorem 3.1.} {\it Under the Basic Assumption 3.0, (i) If n$\geq 2$, then
the following identity holds,
$$\chi(M)=-\sum_{p\in Z_K}{\rm deg}\left(v_K(p)|_{\partial
B_p(\epsilon)}\right).\eqno(3.4)$$ (ii) If $n=1$, set $Z_{K,+}=\{
x\in Z_K: \xi(x),\ \eta(x)\ {\rm form \ an\ oriented \ frame\ at}\
x\}$, then
$$\chi(M)=-\sum_{p\in Z_K\setminus Z_{K,+}}{\rm deg}\left(v_K(p)|_{\partial
B_p(\epsilon)}\right).\eqno(3.4)'$$}

\noindent{\bf Proof.} For simplicity, we first assume that $M$ is
spin and denote by $S(TM)=S_+(TM)\oplus S_-(TM)$ the bundle of
spinors associated with $(TM,g^{TM})$.

Let $v_K=\tau c(\xi)+\sqrt{-1}c(\eta):S_+(TM)\rightarrow S_-(TM)$ be
defined similarly as in (3.3), only that now it is defined on the
whole manifold $M$.

Let $Z({v_K})$ denote   the set   of points at which $v_K$ is not
invertible.

%We first establish the following key lemma.

$\ $

\noindent {\bf Lemma 3.2.} {\it One has, (i) If $n\geq 2$, then $Z(v_K)=Z_K$; (ii) If $n=1$, then $Z(v_K)=Z_K\setminus Z_{K,+}$.}

\noindent{\bf Proof.} From (3.1) and (3.2), it is clear that $p\in
Z_K$ if and only if $|\xi|=|\eta|$ and $\langle \xi,\eta\rangle=0$.

Let $v^*_K:S_-(TM)\rightarrow S_+(TM)$ be the adjoint of $v_K$ with
respect to the natural Hermitian metrics on $S_\pm(TM)$. Set
$V_K=v_K+v^*_K:S(TM)\rightarrow S(TM)$. Then $v_K$ is not invertible
if and only if $V^2_K$ is not strictly positive.

Clearly,
$$V_K=\tau
c(\xi)+\sqrt{-1}c(\eta):S(TM)\rightarrow S(TM).\eqno(3.5)$$ From
(3.5), one finds
$$V^2_K=|\xi|^2+|\eta|^2+\sqrt{-1}\tau
(c(\xi)c(\eta)-c(\eta)c(\xi)).\eqno(3.6)$$

Now if at some  $x\in M$, $|\xi |=|\eta |$ and $\langle
\xi,\eta\rangle=0$, then $V_K^2=2|\xi|^2+2\sqrt{-1}\tau
c(\xi)c(\eta) $ which is clearly seen not invertible if $n\geq 2$ or if $n=1$ but
$\xi$ and $\eta$ do not form an oriented frame at $x$.\footnote{As
one verifies in this case   that either $\xi=\eta=0$, or $c(\xi)-\sqrt{-1}\tau c(\eta)\neq 0$ while $(|\xi|^2+\sqrt{-1}\tau
c(\xi)c(\eta))(c(\xi)-\sqrt{-1}\tau c(\eta))=0$.}

Thus, one gets $Z_K\setminus Z_{K,+}\subset Z(v_K)$.

On the other hand,  observe that if $|\xi|\neq |\eta|$, then
$|\xi|^2+|\eta|^2>2|\xi|\cdot |\eta|$, while it is clear that
$2|\xi|\cdot |\eta|+\sqrt{-1}\tau
(c(\xi)c(\eta)-c(\eta)c(\xi))\geq 0$.

Thus if $|\xi(x)|\neq |\eta(x)|$, then $x$ is not in $Z(v_K)$.

Now if at some $x\in M$, $|\xi |=|\eta |$ and $\langle
\xi,\eta\rangle\neq 0$, one has
$$c(\xi)c(\eta)-c(\eta)c(\xi)=c(\xi)c\left(\eta-{\langle\eta,\xi\rangle\over
|\xi|^2}\xi\right)-c\left(\eta-{\langle\eta,\xi\rangle\over
|\xi|^2}\xi\right)c(\xi),\eqno(3.7)$$ with
$$\left|\eta-{\langle\eta,\xi\rangle\over
|\xi|^2}\xi\right|<|\eta|.\eqno(3.8)$$

From (3.6)-(3.8), one finds that if at some $x\in M$, $|\xi |=|\eta
|$ and $\langle \xi,\eta\rangle\neq 0$, then $V_K^2>0$.

Thus, $Z(v_K)\subset Z_K$. Moreover, if $n=1$, then one verifies
directly that $Z(v_K)\subset Z_K\setminus Z_{K,+}$. The proof of
Lemma 3.2 is completed.\ \ Q.E.D.

$\ $

Back to the proof of Theorem 3.1. By Lemma 3.2, we know that the
Basic Assumption 0.1 holds for $v_K: S_+(TM)\rightarrow S_-(TM)$.
Thus one may apply Theorem 1.1 to it to get
$$\left\langle {\rm ch}\left(S_+(TM)\right)-{\rm
ch}\left(S_-(TM)\right),[M]\right\rangle =(-1)^{n-1}\sum_{p\in
Z(v_K)}{\rm deg}\left(v_K(p)|_{\partial
B_p(\epsilon)}\right).\eqno(3.9)$$

On the other hand, it is standard that (cf. \cite{LM})
$$\left\langle {\rm ch}\left(S_+(TM)\right)-{\rm
ch}\left(S_-(TM)\right),[M]\right\rangle=(-1)^n
\chi(M).\eqno(3.10)$$

From (3.9) and (3.10), one gets (3.4).

Thus we have proved Theorem 3.1 in the case where $M$ is spin.

For the general case where $M$ need not be spin, we may consider the
Signature complex (cf. \cite{LM}) associated with $(TM,g^{TM})$
instead. Then the same argument  above leads to formulas similar  to
(3.9) and (3.10), with the   right hand sides both be multiplied by
a factor $2^n$, while in the left hand sides the Spin complex be
replaced by the Signature complex. Thus one gets again  (3.4). We
leave the details to the interested reader.

The proof of Theorem 3.1 is completed. \  \ Q.E.D.

$\ $

\noindent{\bf Remark 3.3.} If $Z_K=\emptyset$, then one recovers
(and at the same time gives a new proof of) the vanishing result of
Jacobowitz \cite{J} which has been stated in Proposition 1.2.

$\ $

\noindent{\bf Remark 3.4.} Theorem 3.1, in its most general form,
should be regarded as a geometric result. As a simple amazing
consequence (actually a consequence of Proposition 1.2), if
$\chi(M)\neq 0$ and $K=\xi+\sqrt{-1}\eta\in\Gamma(T_{\bf C}M)$ is
{\it nowhere zero} over $M$, then for any Riemannian metric $g^{TM}$
on $TM$, there is at least one point $x\in M$, at which one has
$|\xi|_{g^{TM}}=|\eta|_{g^{TM}}$ and
$\langle\xi,\eta\rangle_{g^{TM}}=0$. Moreover, if $n=1$, then there
exists at least {\it two} such
 points.\footnote{This is because one can switch $\xi$ and $\eta$.}

$\ $

\noindent{\bf Remark 3.5.} One may also extend Theorem 3.1 to the
case where $TM$ is replaced by an arbitrary  oriented Euclidean
vector bundle. We leave the details to the interested reader.

$\ $

Next, we show that Theorem 3.1 is indeed a generalization of the
original Poincar\'e-Hopf index formula (cf. \cite[Theorem
11.25]{BT}).

To do so, we take $\xi=0$, then $Z_K$ is the zero set of $\eta$,
which we have assumed to consist of isolated points.

Without loss of generality we also assume that $|\eta|=1$ on each
$\partial B_p(\epsilon)$, $p\in Z_K$.

In view of the last equality in the proof of Lemma 2.2, one has
$${\rm deg}\left(v_K(p)|_{\partial B_p(\epsilon)}\right)=
{1\over{(2\pi{\sqrt{-1}})^n}}{{
(n-1)!}\over{(2n-1)!}}\int_{\partial{B_{p}(\epsilon)}}{\rm
tr}_{S_+(TM)}\left[\left(v^{-1}dv\right)^{2n-1}\right],\eqno(3.11)$$
with
$$v=\sqrt{-1}c\left(\eta|_{\partial
B_p({\epsilon})}\right).\eqno(3.12)$$

Let $f_1,\cdots,f_{2n-1}$ be an orthonormal basis of $T(\partial
B_p(\epsilon))$, let $f_1^*,\cdots,f_{2n-1}^*$ be the metric dual
basis of $T^*(\partial B_p(\epsilon))$.

From  (3.12), one deduces that (compare with \cite[(27)]{Z})
$${\rm
tr}_{S_+(TM)}\left[\left(v^{-1}dv\right)^{2n-1}\right]=-{2^{n-1}(2n-1)!
(\sqrt{-1})^n}f_1^*\wedge\cdots\wedge f^*_{2n-1}\int^B
\eta^*\wedge\left(\nabla_{f_1}^{TM}\eta\right)^*\wedge\cdots\wedge\left(\nabla_{f_{2n-1}}^{TM}\eta\right)^*,\eqno(3.13)$$
where $\nabla^{TM}$ is the Levi-Civita connection of $g^{TM}$ and
where $\int^B\eta^*\wedge (\nabla_{f_1}^{TM}\eta
)^*\wedge\cdots\wedge (\nabla_{f_{2n-1}}^{TM}\eta )^*$ is the
function on $\partial B_p(\epsilon)$ such that
$$\eta^*\wedge\left(\nabla_{f_1}^{TM}\eta\right)^*\wedge\cdots\wedge\left(\nabla_{f_{2n-1}}^{TM}\eta\right)^*
=\left(d{\rm vol}_{g^{TM}}\right)\int^B
\eta^*\wedge\left(\nabla_{f_1}^{TM}\eta
^*\right)\wedge\cdots\wedge\left(\nabla_{f_{2n-1}}^{TM}\eta\right)^*\eqno(3.14)$$
on $\Lambda^{2n}(T^*M)|_{\partial B_p(\epsilon)}$.

Let $\eta_p:\partial B_p(\epsilon)\rightarrow S^{2n-1}(1)$ denote
the canonical map induced by $\eta|_{\partial B_p(\epsilon)}$.

 By (3.14), one
finds
$$f_1^*\wedge\cdots\wedge f^*_{2n-1}\int^B
\eta^*\wedge\left(\nabla_{f_1}^{TM}\eta\right)^*\wedge\cdots\wedge\left(\nabla_{f_{2n-1}}^{TM}\eta\right)^*
=\eta_p^*\,\omega,\eqno(3.15)$$ where $\omega$ is the volume form on
$S^{2n-1}(1)$.

From (3.11), (3.13) and (3.15), one   gets
$${\rm deg}\left(v_K(p)|_{\partial B_p(\epsilon)}\right)=-{1\over{(2\pi{\sqrt{-1}})^n}}{{
(n-1)!}\over{(2n-1)!}}{2^{n-1}(2n-1)!  (\sqrt{-1})^n}\int_{\partial
B_p(\epsilon)}\eta_p^*\,\omega$$
$$={-(n-1)!\over 2\pi^n}\int_{\partial
B_p(\epsilon)}\eta_p^*\,\omega=-{\rm
deg}\left(\eta_p\right),\eqno(3.16)$$ where ${\rm deg} (\eta_p )$
denotes the Brouwer degree (cf. \cite{BT}) of the map
$\eta_p:\partial B_p(\epsilon)\rightarrow S^{2n-1}(1)$.

From (3.4) and (3.16), one gets
$$\chi(M)=\sum_{p\in {\rm zero\ set\ of}\ \eta}{\rm deg}
\left(\eta_p\right),$$ which is exactly the original Poincar\'e-Hopf
index formula (cf. \cite[Theorem 11.25]{BT}).

$\ $

\noindent{\bf Remark 3.6.} Continuing Remark 3.4 and assume $n\geq
2$. Let $K=\xi+\sqrt{-1}\eta$ be such that the zero set of $\xi$ is
discrete and that $p\in M$ is a zero point of $\xi$ such that ${\rm
deg} (\xi_p )\neq \chi(M)$, while $\eta$ vanishes on a closed ball
of a sufficiently small positive radius around $p$ and is nowhere
zero outside  this closed ball.\footnote{The existence of such a
vector field is clear, as according to a famous theorem of Hopf,
there always exists a vector field on $M$ which vanishes only at
$p$.} Then according to (3.16), $-{\rm deg}(v_K(p))={\rm deg} (\xi_p
)\neq \chi(M)$. Combining this with Theorem 3.1, we see that for any
Riemannian metric $g^{TM}$, there is $x\in M$ such that
$|\xi|_{g^{TM}}=|\eta|_{g^{TM}}\neq 0$ and
$\langle\xi,\eta\rangle_{g^{TM}}=0$. This extends Remark 3.4 to the
case where $K=\xi+\sqrt{-1}\eta$ might vanish on $M$.

$\ $

Now we exhibit an example to illustrate the last line in Remark 3.4.

$\ $

\noindent{\bf Example 3.7.} Let $S^2= \{(x,y,z):x^2+y^2+z^2=1\}$ be
the standard two sphere in the Euclidean space ${\bf R}^3$. Set
$\xi=(-y,x,0)$ and $\eta=(z,0, -x)$. Clearly, as $x^2+y^2+z^2=1$,
$\xi+\sqrt{-1}\eta$ is nowhere zero on $S^2$. Now $|\xi|=|\eta|$
together with $\langle\xi,\eta\rangle=0$ imply that  $x=\pm 1,\
y=z=0$. Thus, $Z_K$ consists of {\it two} points   $p=(1,0,0)$,
$q=(-1,0,0)$. One then verifies that at $q\in S^2$, $\xi=(0,-1,0)$
and $\eta=(0,0,1)$ form an oriented frame of $T_pS^2$. Thus, by
(3.4), one sees that the degree at $p$ equals to $-2$, as the Euler
number of $S^2$  is 2.

 $\ $

 Finally, with the help of Example 3.7, we exhibit an application of Theorem 1.1 in the higher dimensional case.

 $\ $

 \noindent{\bf Example 3.8.} We take a product $M=S^2\times \cdots\times S^2$ with $m\geq 2$ copies of $S^2$.
 We use a subscript to denote the corresponding factor of $S^2$. So now let $\xi_i,\ \eta_i$, $1\leq i\leq m$,
 be the vector fields constructed in Example 3.7 on the $i$-th factor $S^2$ (denoted by $S^2_i$). Let $v_{K,i}$
 be the lifting to $M$ of the corresponding map defined as in the proof of Theorem 3.1 on
 $S^2_i$. Then each $v_{K,i}$ maps $\Gamma(S_+(TM))$ to
 $\Gamma(S_-(TM))$.
 Set $v_K=\sum_{i=1}^m v_{K,i}$, then one verifies directly that $v_K$ is singular only at the point
 $(p_1,\cdots,p_m)\in S^2\times \cdots\times S^2$. By combining Theorem 1.1 with (3.9) and (3.10), one then
 gets that the degree of $v_K$ at $(p_1,\cdots,p_m)$ equals to $-2^m$, as the Euler number
 of $S^2\times \cdots\times S^2$ equals to $2^m$. Conversely, one can compute   the
 degree at $(p_1,\cdots,p_m)$ first,  and then get  the Euler number
  of $S^2\times \cdots\times S^2$ by using Theorem 1.1.

$$\ $$

\noindent H. Feng, School of Mathematics and Statistics, Chongqing
University of Technology, Chongqing 400050, PR China

\noindent{\it Email}: fht@nankai.edu.cn

$\ $

\noindent W. Li, Department of Mathematics, Oklahoma State
University, Stillwater, OK 74078-1058, USA

\noindent{\it Email}: wli@math.okstate.edu

$\ $

\noindent W. Zhang, Chern Institute of Mathematics \& LPMC, Nankai
University, Tianjin 300071, PR China

\noindent{\it Email}: weiping@nankai.edu.cn

\end{document}